\documentclass[a4paper,12pt]{article}
\usepackage[T1]{fontenc}
\usepackage[english]{babel}
\usepackage[dvips]{graphicx}
\usepackage{amssymb,amsthm,amsmath}
\usepackage{pstricks}

\newcommand{\atan}{\textnormal{atan}\,}

\newcommand{\card}{\textnormal{card}\,}
\newcommand{\e}{\varepsilon}
\newcommand{\J}{ {B_j} }
\newcommand{\PP}{ { \R^2 \setminus\{ 0 \} } }
\newcommand{\PS}{ { \Rn \setminus\{ 0 \} } }
\newcommand{\R}{\mathbb{R}}
\newcommand{\C}{\mathbb{C}}
\newcommand{\Rn}{ {\mathbb{R}^n} }

\newcommand{\cmmnt}[1]{}

\newtheorem{theorem}[equation]{Theorem}
\newtheorem{lemma}[equation]{Lemma}
\newtheorem{proposition}[equation]{Proposition}
\newtheorem{corollary}[equation]{Corollary}

\theoremstyle{definition}
\newtheorem{remark}[equation]{Remark}
\newtheorem{definition}[equation]{Definition}
\newtheorem{example}[equation]{Example}
\newtheorem{questions}[equation]{Questions}
\numberwithin{equation}{section}

\newcounter{minutes}
\divide\time by 60
\newcounter{hours}
\multiply\time by 60 \addtocounter{minutes}{-\time}
\begin{document}
\psset{linewidth=1pt}

\begin{center}
{\Large \bf Local Convexity Properties of $j$-metric Balls}
\end{center}
\medskip

\begin{center}
{\large Riku Klén }
\end{center}
\bigskip

\begin{abstract}
  This paper deals with local convexity properties of the $j$-metric. We consider convexity and starlikeness of the $j$-metric balls in convex, starlike and general subdomains of $\Rn$.

  \hspace{5mm}


\end{abstract}
\bigskip
\begin{center}
\texttt{File:~\jobname .tex, 2007-01-10,
         printed: \number\year-\number\month-\number\day,
         \thehours.\ifnum\theminutes<10{0}\fi\theminutes}
\end{center}

%%%%%%%%%%%%%%%%%%%%%%%%%%%%%%%%%%%%%%%%%%%%%%%%%%%%%%%%%%%%%
\section{Introduction}

The \emph{$j$-distance} in a proper subdomain $G$ of the Euclidean space $\Rn$, $n \ge 2$, is defined by
$$
  j_G(x,y) = \log \left( 1+\frac{|x-y|}{\min \{ d(x),d(y) \}}\right),
$$
where $d(x)$ is the Euclidean distance between $x$ and $\partial G$. If the domain $G$ is understood from the context we use notation $j$ instead of $j_G$.

The $j$-distance was first introduced by F.W. Gehring and B.P. Palka \cite{gp} in 1976 in a slightly different form and in the above form, by M. Vuorinen \cite{vu2} in 1985. The $j$-distance is actually a metric and a proof of the triangle inequality valid for general metric spaces is given in \cite{s}. Previously the $j$-metric has been studied in connection with the study of other metrics \cite{go,h,s,v,vu2}. See also recent papers \cite{hl,l}. In spite of these studies many basic questions of the $j$-metric remain open and some of them will be studied here.

The purpose of this paper is to study metric spaces $(G,j_G)$ and especially local convexity properties of \emph{$j$-metric balls} or in short \emph{$j$-balls} defined by
$$
  \J(x,M) = \{ y \in G \colon j(x,y) < M \},
$$
where $M > 0$ and $x \in G$. In the dimension $n=2$ we call these \emph{$j$-metric disks} or \emph{$j$-disks}.

M. Vuorinen suggested in \cite{vu4} a general question about the convexity of balls of small radii in metric spaces. This paper is motivated by this question and we will provide an answer in a particular case. Our main result is the following theorem. For the definition of starlike domains see \ref{defstarlikeness}.

\begin{theorem}\label{convexj}
  For a domain $G \subsetneq \Rn$ and $x \in G$ the $j$-balls $\J(x,M)$ are convex if $M \in (0,\log 2]$ and strictly starlike with respect to $x$ if $M \in \big( 0,\log (1+\sqrt 2) \big]$.
\end{theorem}

In Section \ref{Properties} we consider general properties of the $j$-metric and show that for any $G$ there exists points such that there is no geodesic between them. In Section \ref{PuncturedS} we consider local convexity properties of $j$-balls in punctured space and in Section \ref{GeneralG} we extend these results to an arbitrary domain $G \subsetneq \Rn$. We will further consider convexity of $j$-balls in convex domains and starlikeness of $j$-balls in starlike domains.

%%%%%%%%%%%%%%%%%%%%%%%%%%%%%%%%%%%%%%%%%%%%%%%%%%%%%%%%%%%%%
\section{Properties of the $j$-metric}\label{Properties}

Throughout this paper $G \subsetneq \Rn$, $n \ge 2$, is a domain. We denote $m(x,y) = \min \{ d(x),d(y) \}$ and we use notation $B^n(x,M)$ for the Euclidean balls and $S^{n-1}(x,M)$ for the Euclidean spheres. We often identify $\R^2$ with the complex plane $\C$.

In 1976 F.W. Gehring and B.P. Palka \cite{gp} also introduced the quasihyperbolic metric, which has been widely applied in geometric function theory and mathematical analysis in general, see e.g. \cite{vu3,v}. The \emph{quasihyperbolic distance} between two points $x$ and $y$ in a proper subdomain $G$ of the Euclidean space $\Rn$, $n \ge 2$, is defined by
$$
  k_G(x,y) = \inf_{\alpha \in \Gamma_{xy}} \int_{\alpha}\frac{|dx|}{d(x)},
$$
where $\Gamma_{xy}$ is the collection of all rectifiable curves in $G$ joining $x$ and $y$. We denote the \emph{quasihyperbolic ball} by
$$
  D_G(x,M) = \{ y \in G \colon k_G(x,y) < M \}.
$$

The quasihyperbolic metric is closely related with the $j$-metric. By \cite[Lemma 2.1]{gp} $j_G$ is always a minorant of $k_G$, in other words, for a proper subdomain $G$ of $\Rn$ we have
$$
  j_G(x,y) \le k_G(x,y)
$$
for all $x,y\in G$.

The following result can be used to estimate the quasihyperbolic
metric from above by the $j$-metric.
\begin{proposition}\cite[Lemma 3.7]{vu3}
  Let $G \subsetneq \Rn$ be a domain, $x \in G$, $y \in B^n \big( x,d(x) \big)$ and $s \in (0,1)$. Then
  $$
    k_G(x,y) \le \frac{1}{1-s} j_G(x,y).
  $$
\end{proposition}

The following lemma gives Euclidean bounds for the $j$-balls.

\begin{proposition}\label{roundness}\cite[Theorem 3.8]{s}
  For a proper subdomain $G \subset \Rn$, $x \in G$ and $M > 0$ we have
  $$
    B^n \big( x,r \, d(x) \big) \subset B_j(x,M) \subset B^n \big( x,R \, d(x) \big),
  $$
  where $r =1- e^{-M}$ and $R = e^M-1$. Moreover
  $$
    \J(x,M) \subset \left\{ z \in G \colon e^{-M}d(x) \le d(z) \le e^M d(x) \right\}.
  $$
\end{proposition}

\begin{remark}
  A similar result to Proposition \ref{roundness} is also true for the quasihyperbolic metric see \cite[page 347]{vu1}.
\end{remark}

By Proposition \ref{roundness} the $j$-ball $\J(x,M)$ shrinks towards the center $\{ x \}$ as $M$ approaches 0. The following lemma shows that the $j$-balls $\J(x,M)$ exhaust the domain $G$.

\begin{lemma}
  Let $G \subset \Rn$ be a bounded domain and fix $x \in G$ and $s \in (0,d(x)]$. Then
  $$
    \{ y \in G \colon d(y) > s \} \subset \J \big( x,\log(1+d/s) \big),
  $$
  for $\displaystyle d = \sup_{z \in \partial G} |x-z|$.
\end{lemma}
\begin{proof}
  Let us assume $d(y) > s$. Then either $m(x,y) = d(x) \ge s$ or $m(x,y) = d(y) > s$. In both cases $m(x,y) \ge s$ and since $|x-y| < d$ for all $y \in G$ we have
  $$
    j(x,y) = \log \left( 1+\frac{|x-y|}{m(x,y)} \right) < \log \left( 1+\frac{d}{s} \right).
  $$
\end{proof}

Let us denote the set of closest boundary points of a point $x$ in a domain $G \subset \Rn$ by
$$
  R_x = \left\{ z \in \partial G \colon |z-x| = d(x) \right\}.
$$

The next result characterizes the case of equality in the triangle inequality for the $j$-metric. Its proof is based on the proof of the triangle inequality \cite[Lemma 2.2]{s}.

\begin{theorem}\label{trianglej}
  Let $x,y,z \in G \subsetneq \Rn$ be distinct points and $d(x) \le d(z)$. Then
  $$
    j_G(x,z) = j_G(x,y)+j_G(y,z)
  $$
  implies that $x$, $z$ and $u$ are collinear for some $u \in R_x$ and $y \in (x,z)$ with $d(x) < d(y) < d(z)$.
\end{theorem}
\begin{proof}
  By definition $j_G(x,z) < j_G(x,y)+j_G(y,z)$ is equivalent to
  \begin{equation}\label{trianglein}
    \frac{|x-z|}{m(x,z)} < \frac{|x-y|}{m(x,y)}+\frac{|y-z|}{m(y,z)}+\frac{|x-y| |y-z|}{m(x,y) m(y,z)}.
  \end{equation}
  The assumption $d(x) \le d(z)$ implies $m(x,z) = d(x)$.

  If $d(y) \le d(x)$, then the inequality (\ref{trianglein}) is equivalent to
  $$
    |x-z| < |x-y|\frac{d(x)}{d(y)} + |y-z|\frac{d(x)}{d(y)} + \frac{|x-y||y-z|}{d(y)} \frac{d(x)}{d(y)},
  $$
  which is true, because $|x-z| \le |x-y|+|y-z|$, $(|x-y||y-z|)/d(y) > 0$ and $d(x)/d(y) \ge 1$.

  If $d(y) > d(x)$, then the inequality (\ref{trianglein}) is equivalent to
  $$
    |x-z| < |x-y| + |y-z| \left( \frac{d(x)+|x-y|}{m(y,z)} \right),
  $$
  which is false if and only if $x$, $y$ and $z$ are collinear and
  $$
    \frac{d(x)+|x-y|}{m(y,z)} = 1.
  $$
  If $d(x) = d(z)$, then $d(x)/m(y,z) = 1$ and

  \begin{equation}\label{ineq}
    \frac{d(x)+|x-y|}{m(y,z)} > 1.
  \end{equation}
  If $d(x) < d(z) < d(y)$, then the inequality $(\ref{ineq})$ is true, because $d(x) + |x-y| \ge d(y) > d(z) = m(y,z)$. If $d(x) < d(y) \le d(z)$, then the inequality $(\ref{ineq})$ is true if and only if $y \notin \{ k(x-u) \colon k > 0 \}$, where $u \in R_x$.
\end{proof}

The implication of Theorem \ref{trianglej} in the other direction was proved by Hästö, Ibragimov and Lindén \cite[Corollary 3.7]{hil}.

\begin{definition}
  Let $G \subsetneq \Rn$ be a domain and $\gamma$ a curve in $G$. If
  $$
    j(x,y) + j(y,z) = j(x,z)
  $$
  for all $x,z \in \gamma$ and $y \in \gamma'$, where $\gamma'$ is the subcurve of $\gamma$ joining $x$ and $z$, then $\gamma$ is a \emph{geodesic segment} or shortly a \emph{geodesic}. We denote a geodesic between $x$ and $y$ by $J[x,y]$.
\end{definition}

By Theorem \ref{trianglej} and the result of Hästö, Ibragimov and Lindén we can easily find all geodesics $J[x,y]$ for any domain $G$. The geodesic needs to satisfy the triangle inequality as equality at each point and therefore the geodesic can only be a line segment $l$ with the following property.

\begin{lemma}\label{geodesiclemma}
  Let $G \subsetneq \Rn$ be a domain and $J[x,y]$ be a geodesic segment with $x,y \in G$. There exists $u \in \partial G$ such that $u \in R_s$ for all $s \in J[x,y]$ and $u$ and $J[x,y]$ are collinear.
\end{lemma}
\begin{proof}
  Let us assume, on the contrary, that there exists $z \in J[x,y]$ such that $d(z) < d(x)-|x-z|$. Now $j_G(x,z) + j_G(z,y) = j_G(x,y)$ is equivalent to
  $$
    d(z)|x-z| + \big( d(x)+|x-z| \big)|z-y| = d(z)|x-y|.
  $$
  We have
  \begin{eqnarray*}
    d(z)|x-y| & \le & d(z)|x-z| + d(z)|z-y|\\
    & < & d(z)|x-z| + \big( d(x)+|x-z| \big)|z-y|\\
    & = & d(z)|x-y|
  \end{eqnarray*}
  which is a contradiction.
\end{proof}

\begin{theorem}\label{notgeodesic}
  Let $G \subsetneq \Rn$ be a domain. Then there exist $x,y \in G$ such that there is no geodesic $J[x,y]$.
\end{theorem}
\begin{proof}
  Let us assume, on the contrary, that for all $x,y \in G$ there exists a geodesic $J[x,y]$. Since $G$ is a domain, we can choose $x,y,z \in G$ to be three distinct noncollinear points. Now there exists a geodesic $J[x,y]$ from $x$ to $y$. We may assume $d(x) < d(y)$ and then by Lemma \ref{geodesiclemma} $B^n \big( x,d(x) \big) \subset B^n \big( y,d(y) \big) \subset G$.

  On the other hand, there exists a geodesic $J[x,z]$ from $x$ to $z$ and therefore there has to exist a point $u \in S^{n-1} \big( x,d(x) \big) \cap \partial G$ such that $x$, $z$ and $u$ are collinear. This is a contradiction, because $x$, $y$ and $u$ are noncollinear and therefore $u \in B^n \big( y,d(y) \big)$.
\end{proof}

\begin{remark}
  By Theorem \ref{notgeodesic} a $j$-metric geodesic does not always exist between two points. F.W. Gehring and B.G. Osgood have proved \cite[Lemma 1]{go} that for the quasihyperbolic metric there always exists a geodesic between two points of a domain $G \subsetneq \Rn$. 

  However, the geodesics of the $j$-metric are unique while the geodesics of the quasihyperbolic metric need not be unique.
\end{remark}

%%%%%%%%%%%%%%%%%%%%%%%%%%%%%%%%%%%%%%%%%%%%%%%%%%%%%%%%%%%%%
\section{Convexity and starlikeness of $j$-balls in punctured space}\label{PuncturedS}

In this section we consider the case $G = \PS$. By definition the $j$-balls in punctured space $G = \PS$ are similar, which means that $\J(x,M)$ can be mapped onto $\J(y,M)$ for all $x,y \in G$ by rotation and stretching. We see easily that these balls are also symmetric along the line that goes through 0 and the center point.

\begin{theorem}\label{punctconvexj}
  Let $x \in \PS$. Then\\
  \noindent1) the $j$-ball $\J(x,M)$ is convex if and only if $M \in (0,\log 2]$.\\
  \noindent 2) the $j$-ball $\J(x,M)$ is strictly convex if and only if $ M \in (0,\log 2)$.
\end{theorem}
\begin{proof}
  1) By similarity we can assume $x = e_1$ and by symmetry it is sufficient to consider only the case $n=2$. We will consider $\partial \J(1,M)$ for fixed $M$. By definition we have for $z \in \partial \J(1,M)$
  $$
    M = \left\{ \begin{array}{l} \log(1+|z-1|),\quad |z| \ge 1,\\ \log \left(1+|z-1|/|z| \right),\quad |z| < 1,  \end{array} \right.
  $$
  which is equivalent to
  $$
    e^M-1 = \left\{ \begin{array}{l} |z-1|,\quad |z| \ge 1,\\ \left| 1-1/z \right|,\quad |z| < 1. \end{array} \right.
  $$
  For $|z| \ge 1$ the $\partial \J(1,M)$ is an arc of a circle with center 1 and radius $e^M-1$. For $|z| < 1$ the $\partial \J(1,M)$ is a circle that goes through points $1/(e^M)$ and $1/(2-e^M)$ and has center on the real axis. This means that the center of the circle is $c = 1/ \big( e^M( 2-e^M) \big)$ and the radius of the circle is $|e^M-1|/|e^M(2-e^M)|$. Now $c > 1$, if $M \le \log 2$, and $c < 0$, if $M > \log 2$. Therefore $\partial \J(1,M)$ is convex for $M \le \log 2$ and not convex for $M > \log 2$.

  2) We have $c \in (1,\infty)$, where $c$ is as above. Therefore $\J(x,M)$ is strictly convex. In the case $M = \log 2$ we have $c = \infty$ and $\J(x,M)$ is not strictly convex.
\end{proof}

\begin{remark}\label{rem1}
  For fixed $x \in G$ the quasihyperbolic ball $D_G(x,M)$ is strictly convex in $G = \PS$ if and only if $M \in (0,1]$ \cite{k}.
\end{remark}

Clearly $\J(x,M)$ is never smooth. We will next define starlikeness of a domain.

\begin{definition}\label{defstarlikeness}
  Let $G \subset \Rn$ be a bounded domain and $x \in G$ . We say that $G$ is \emph{starlike with respect to $x$} if each line segment from $x$ to $y \in G$ is contained in $G$. The domain $G$ is \emph{strictly starlike with respect to $x$} for $x \in G$ if each ray from $x$ meets $\partial G$ at exactly one point.
\end{definition}

The next theorem determines the values of $M$ for which the $j$-ball $\J(x,M)$ is strictly starlike with respect to $x$.

\begin{theorem}\label{starlikedj}
  For $x \in \PS$ the $j$-ball $\J(x,M)$ is strictly starlike with respect to $x$ if and only if $M \in \big( 0, \log (1+\sqrt{2}) \big]$.
\end{theorem}
\begin{proof}
  Because the $j$-balls are similar it is sufficient to consider $x = e_1$. By symmetry it is sufficient to consider the case $n=2$ and the part of $\partial \J(1,M)$ that is above the real axis. If $M \ge \log 3$, then $\J(1,M) = B^2(1,r) \setminus B^2(c,s)$, where $c$, $r$ and $s$ are given in the proof of Theorem \ref{punctconvexj} and $B^2(c,s) \subset B^2(1,r)$. Therefore $\J(1,M)$ can be starlike with respect to 1 only for $M < \log 3$.

  Let us assume $M < \log 3$. By the proof of Theorem \ref{punctconvexj} $\J(1,M) = B^2(1,r) \setminus B^2(c,s)$. Let us denote the point of intersection of $S^1(1,r)$ and $S^1(c,s)$ above the real axis by $z$. Now $z$ is also the point of intersection of the unit circle and the boundary $\partial \J(1,M)$. Let us denote by $l$ the line that goes through points 1 and $z$. Now $\J(1,M)$ is strictly starlike with respect to 1 if and only if $l \cap B^2(1,r) \cap B^2(c,s) = \emptyset$. If $z$ is a tangent of $S^1(c,s)$, then the circles $S^1(1,r)$ and $S^1(c,s)$ are perpendicular and $M$ has the largest value such that $\J(1,M)$ is starlike with respect to 1.

  By the proof of Theorem \ref{punctconvexj} we have $c = -1/e^M(e^M-2)$, $r = |1-z| = e^M-1$, $|1-c| = (e^M-1)^2/e^M(e^M-2)$ and $s = |z-c| = (e^M-1)/e^M(e^M-2)$. Let us assume that $z$ is a tangent of $S^1(c,s)$. Now by the Pythagorean Theorem
  $$
    \frac{(e^M-1)^4}{e^{2M}(e^M-2)^2} = (e^M-1)^2+\frac{(e^M-1)^2}{e^{2M}(e^M-2)^2},
  $$
  which is equivalent to $e^{2M}-2e^M-1=0$ and therefore
  $$
    M = \log (1+\sqrt{2}).
  $$
\end{proof}

\begin{figure}[htp]
  \begin{center}
    \includegraphics[width=5cm]{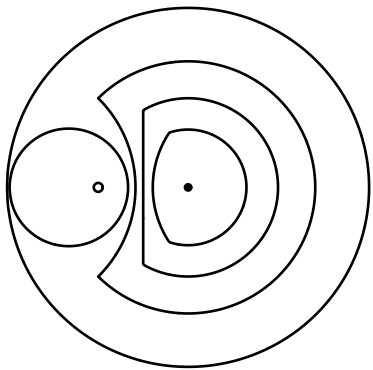}
    \caption{The boundaries of $j$-disks $j(1,M)$ in punctured plane $G = \PP$ with $M=0.5$, $M=\log 2$, $M=\log(1+\sqrt{2})$ and $M=1.1 \approx \log 3$.}
  \end{center}
\end{figure}

\begin{example}
  Let us consider the starlikeness of $j$-balls $\J(x,M)$ with respect to $z \in \J(x,M)$ for $M > \log 2$. By choosing $z = (e^{-M}+\e)x/|x|$ for $\e > 0$ and letting $\e$ approach to zero we see that $\J(x,M)$ is not starlike with respect to $z$.

  On the other hand, if we choose $z = (e^M-\e)x/|x|$ for $\e > 0$ and $M < \log \big( (3+\sqrt{5}/2 \big)$, we see that $\J(x,M)$ is strictly starlike with respect to $z$ for small enough $\e$.
\end{example}

\begin{remark}\label{rem2}
  For fixed $x \in G$ the quasihyperbolic ball $D_G(x,M)$ is strictly starlike with respect to $x$ in $G = \PS$ if and only if $M \in (0,\kappa]$ \cite{k}, where $\kappa \approx 2.83297$.
\end{remark}

%%%%%%%%%%%%%%%%%%%%%%%%%%%%%%%%%%%%%%%%%%%%%%%%%%%%%%%%%%%%%
\section{Convexity and starlikeness of $j$-balls}\label{GeneralG}

We will consider convexity and starlikeness of $j$-balls $\J(x,M)$ for $M > 0$ in convex, starlike and general domains.

Let us consider $j$-balls in a domain $G$ with a finite number of boundary points. The case $\card \partial G = 1$ is identical to $G = \PS$. If $\partial G = \{ y_1,y_2 \}$, then $B_{j_G}(x,M) = B_{j_{\Rn \setminus \{ y_1 \}}}(x,M) \cap B_{j_{\Rn \setminus \{ y_2 \}}}(x,M)$. This is clear, because the $j$-distance between $a$ and $b$ depends only on the closest boundary point of the end points $a $ and $b$. Similarly for $\partial G = \{ y_1,y_2,\dots ,y_m \}$ we have
$$
    B_{j_G}(x,M) = \bigcap_{i=1}^m B_{j_{\Rn \setminus \{ y_i \}}}(x,M).
$$

\begin{figure}[htp]
  \begin{center}
    \includegraphics[width=27mm]{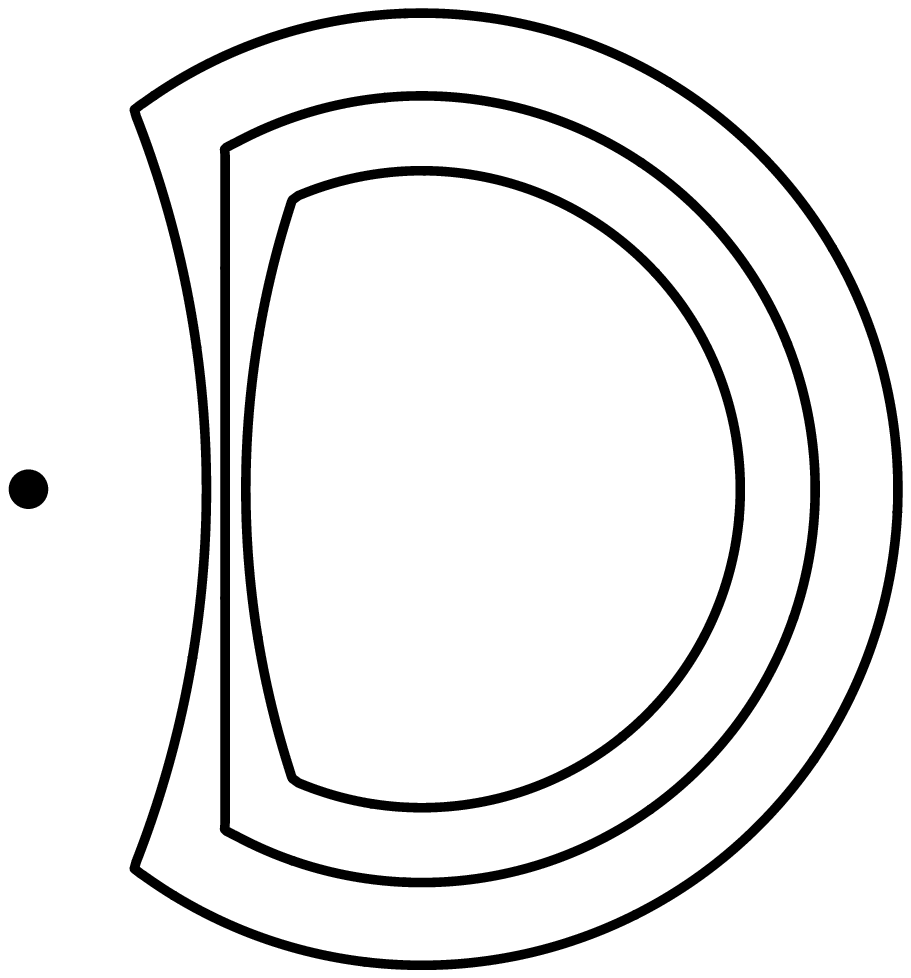}\hspace{6mm}
    \includegraphics[width=27mm]{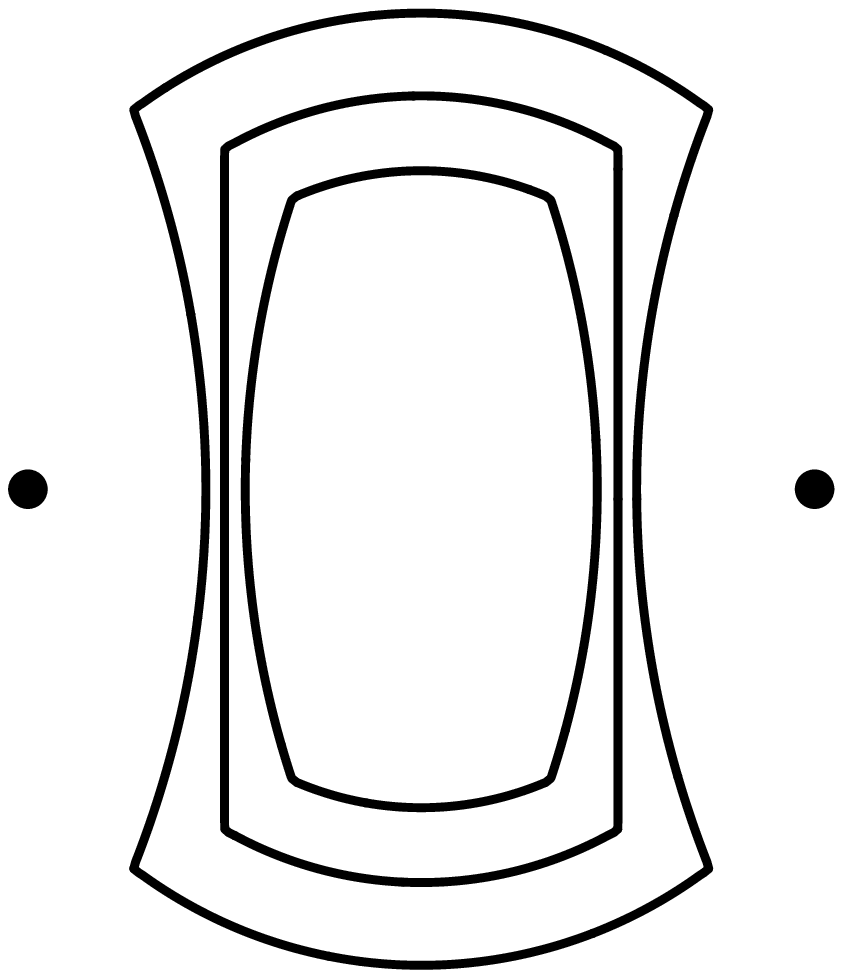}\hspace{6mm}
    \includegraphics[width=27mm]{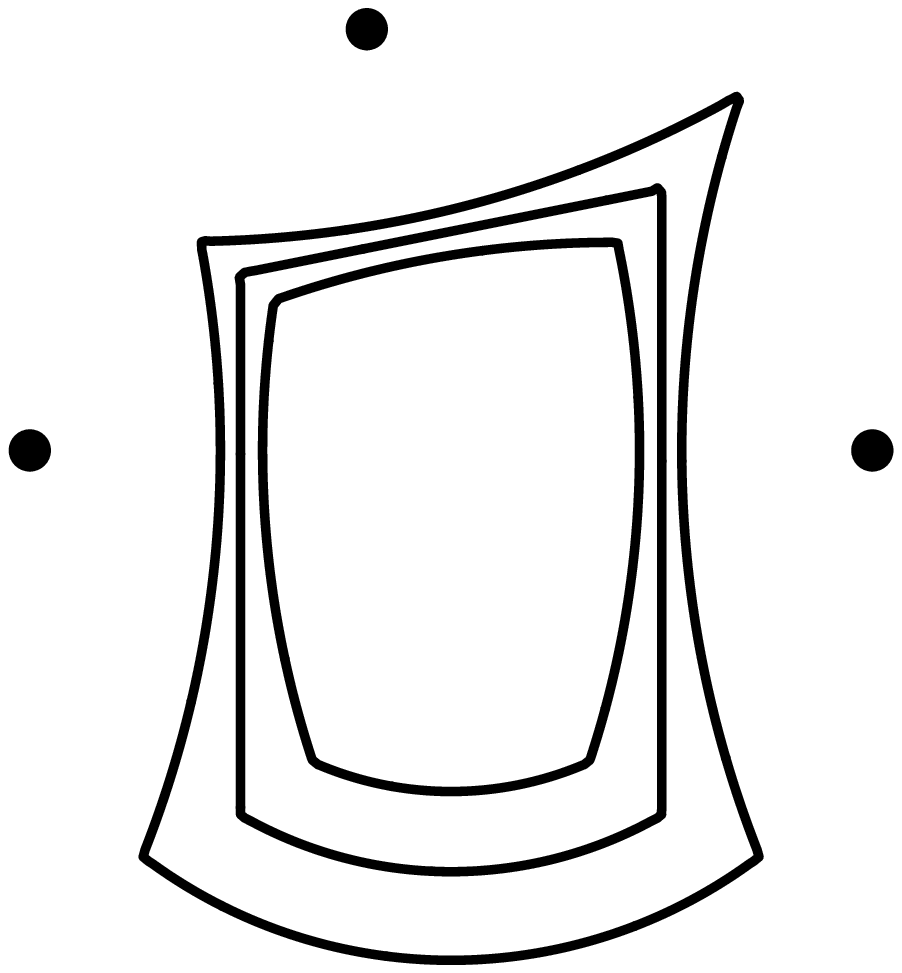}\hspace{6mm}
    \includegraphics[width=27mm]{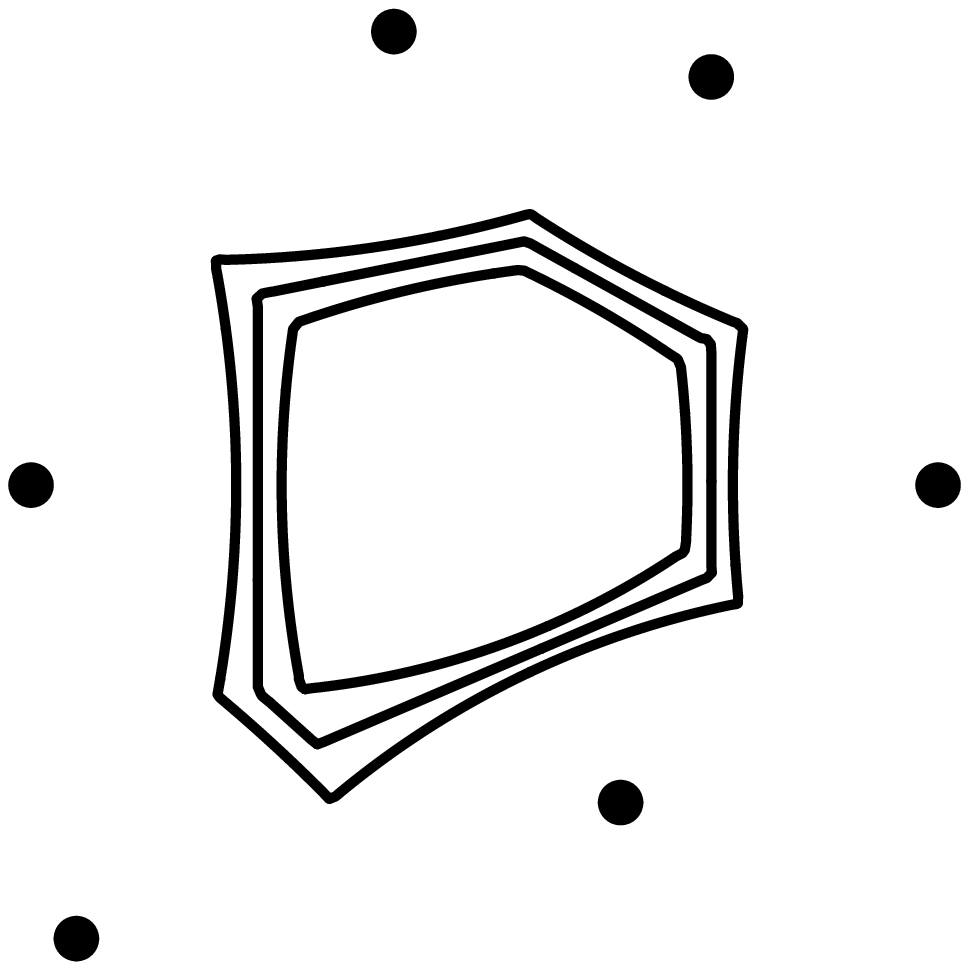}
    \caption{The boundaries of $j$-disks in a domain with 1, 2, 3 and 6 boundary points.}
  \end{center}
\end{figure}

This gives an idea to prove Theorem \ref{convexj}, which shows that $j$-balls are convex in any domain $G$ for small radius $M$.

\begin{proof}[Proof of Theorem \ref{convexj}]
  Let $x \in G$ be arbitrary.  We claim that
  \begin{equation}\label{intersection}
    A = B_{j_G}(x,M) = \bigcap_{z \in \partial G} B_{j_{\Rn \setminus \{ z \}}}(x,M) = B.
  \end{equation}
  Let $y \in B$. We can choose $z' \in \partial G$ with
  $$
    j_{\Rn \setminus \{ z' \}}(x,y) = \min_{z \in \partial G} j_{\Rn \setminus \{ z \}}(x,y).
  $$
  Because $z' \in \partial G$ we have $j_G(x,y) \le j_{\Rn \setminus \{ z' \}}(x,y)$ and therefore $y \in A$.

  On the other hand, let $y \in A$. By definition there is a point $z' \in \partial G$ with $\min \{ |x-z'|,|y-z'| \} = \min_{z \in \partial G} \{ |x-z|,|y-z| \}$. Now $j_{\Rn \setminus \{ z' \}} (x,y) \le j_G(x,y)$ and $y \in B$.

By Theorem \ref{punctconvexj} each $B_{j_{\Rn \setminus \{z\}}}(x,M)$ is convex for $0 < M \le \log 2$ and (\ref{intersection}) $B_{j_G}(x,M)$ is an intersection of convex domains and therefore it is convex.

  If $M \in (0,\log 2]$, then $B_{j_G}(x,M)$ is convex and therefore also starlike with respect to $x$. If $M \in (\log 2,\log (1+\sqrt 2)]$, then
  $$
    \J(x,M) = B \setminus \left( \bigcup_{z \in \partial G} A_z \right),
  $$
  where $B = B^n \big( x,(e^M-1)d(x) \big)$ and $A_z = B^n(c_z z,r_z)$ for $c_z = |z|/(e^M(2-e^M))$ and $r_z = |z||1-e^{-M}|/|e^M-2|$. Let us assume that $\J(x,M)$ is not strictly starlike with respect to $x$. Now there exists $a,b \in B$ such that $b \in (x,a)$, $a \in \J(x,M)$ and $b \notin \J(x,M)$. Now $b \in B^n(c_z z,r_z)$ for some $z \in \partial G$. By the proof of Theorem \ref{starlikedj} $a \in B^n(c_z z,r_z)$, which is a contradiction.
\end{proof}

\begin{corollary}\label{simplycon}
  For a domain $G \subsetneq \Rn$ and $x \in G$ the $j$-balls $\J(x,M)$ are simply connected if $M \in \big( 0,\log (1+\sqrt 2) \big]$.
\end{corollary}
\begin{proof}
  By Theorem \ref{convexj} $B_{j_G}(x,M)$ is starlike with respect to $x$ and therefore also simply connected.
\end{proof}

\begin{corollary}\label{sconvexj}
  For a domain $G \subsetneq \Rn$ and $x \in G$ the $j$-balls $\J(x,M)$ are strictly convex if $M \in (0,\log 2)$.
\end{corollary}
\begin{proof}
  By the proof of Theorem \ref{convexj} and Theorem \ref{punctconvexj}
  $$
    \J(x,M) = \bigcap_{z \in \partial G} (B_{z,1} \cap B_{z,2}),
  $$
  where $B_{z,i}$ is a Euclidean ball and $x \in B_{z,i}$. Therefore $\J(x,M)$ is strictly convex.
\end{proof}

Bounds of Theorem \ref{convexj} are sharp as $G = \PS$ shows. Also the bound $\log (1+\sqrt 2)$ of Corollary \ref{simplycon} is sharp. This can be seen by choosing $G = \R^2 \setminus \{ 0,z \}$ for a certain $z$ and considering $\J(e_1,M)$ for $M > \log (1+\sqrt 2)$. By the proof of Theorem \ref{punctconvexj} we know that
  $$
    \J(e_1,M) = B^2(e_1,r_1) \setminus B^2(c,r_2)
  $$
  for $r_1 = e^M-1$, $c = e_1/ \big( e^M(2-e^M) \big)$ and $r_2 = (e^M-1)/ \big( e^M(e^M-2) \big)$. Let $l$ be the tangent line of $B^2(c,r_2)$ that goes through $e_1$. Denote $\{ y \} = S^1(c,r_2) \cap l$. Choose $z$ to be the reflection of 0 in the line $l$. By a simple computation we have
  $$
    |y-e_1| = \frac{e^M-1}{\sqrt{e^M(e^M-2)}} < r_1.
  $$
  Let us denote by $c'$ the reflection of $c$ in the line $l$. Now $B_{j_{\R^2 \setminus \{ 0,z \}}}(e_1,M) = B^2(e_1,r_1) \setminus \big( B^2(c,r_2) \cup B^2(c',r_2) \big)$ and therefore $\J(e_1,M)$ is disconnected for $M > \log (1+\sqrt{2})$.

  Similar counterexamples can be constructed for $n > 2$. Let us assume $n \ge 2$ and $M > \log (1+\sqrt{2})$. Now we choose
  $$
    G = \Rn \setminus \big( S^{n-1}(z,|z|) \setminus B^n(e_1,1) \big),
  $$
  where $z \in S^{n-1}(e_1,e^M-1)$ and the line $[z,e_1]$ is a tangent of $S^{n-1}(c,r)$ for $c = e_1/ \big( e^M(2-e^M) \big)$ and $r = |1-e^M|/|e^M(2-e^M)|$. Let $y \in [z,e_1] \cap S^{n-1}(e_1,e^M-1)$. Now $j_G(e_1,y) = M$ and $j_G \big( e_1,\frac{1}{2}(z+y) \big) < M$. Therefore $\J(e_1,M)$ is disconnected.

\begin{remark}
  The idea of the proof of Theorem \ref{convexj} cannot be used for the quasihyperbolic metric. We always have
  $$
    D_G(x,M) \subset \bigcap_{z \in \partial G} D_{\Rn \setminus \{ z \}}(x,M)
  $$
  but inclusion in the other direction is not always true. For example $G = \Rn \setminus \{ 0,e_1 \}$, $x = e_1/4$ and $M=1$ gives an counterexample. Now $y = e_1(1-1/e)$ is on the boundary $\partial D_G(x,M)$ because
  $$
    k_G(x,y) = k_{\Rn \setminus \{ 0 \}}(x,e_1/2) + k_{\Rn \setminus \{ e_1 \}}(e_1/2,y) = \log 2 + \log (e/2) = 1.
  $$
  On the other hand, $z = e_1 \big( 1-3/(4 e) \big)$ belongs to the boundary $\partial D_{\Rn \setminus \{ e_1 \}}(x,M)$. Now $0.632 \approx |y| < |z| \approx 0.724$ and therefore $D_{\Rn \setminus \{ 0 \}}(x,M) \cap D_{\Rn \setminus \{ e_1 \}}(x,M) \not \subset D_G(x,M)$.
\end{remark}

The next theorem states convexity of $j$-balls in convex domains.

\begin{theorem}\label{convexdomj}
  Let $M > 0$, $G \subsetneq \Rn$ be a convex domain and $x \in G$. Then $j$-balls $\J(x,M)$ are convex.
\end{theorem}
\begin{proof}
  By Theorem \ref{convexj} we need to consider only the case $M > \log 2$. Let us divide $G$ into two parts $D_1 = \{ z \in G \colon d(z) \ge d(x) \}$ and $D_2= G \setminus D_1$. We will first show that convexity of $G$ implies convexity of $D_1$. Let us assume that $D_1$ is not convex. There exists $a,b \in D_1$ such that $c = (a+b)/2 \notin D_1$ and $d(a) = d(x) =d(b)$. Now $B^n \big( a,d(x) \big)$ and $B^n \big( b,d(x) \big)$ does not contain any points of $\partial G$, but $B^n(c,r)$ for some $r < d(x)$ contains at least one point of $\partial G$. Therefore $G$ is not convex, which is a contradiction.

Let us consider $\J(x,M) \cap D_1$. By definition of the $j$-metric we have for $y \in \partial \J(x,M) \cap D_1$
  $$
    |x-y| = d(x) \big( e^M-1 \big)
  $$
  and therefore $\partial \J(x,M) \cap D_1$ is a subset of $S^{n-1}(x,r)$, where $r=d(x) \big( e^M-1 \big)$. By convexity of $D_1$ the domain $\J(x,M) \cap D_1$ is convex.

  Let us then show that each chord with end points in $\J(x,M) \cap D_2$ is contained in $\J(x,M)$. By definition for $y \in \partial \J(x,M) \cap D_2$ we have
  \begin{equation}\label{ratio}
    d(y) = \frac{|x-y|}{e^M-1}.
  \end{equation}
  Let us assume $y_1,y_2 \in \J(x,M) \cap D_2$ and $z = (y_1+y_2)/2 \notin \J(x,M)$. If $z \in D_1$, then $z \in \J(x,M)$ because $\J(x,M) \subset B^n(x,r)$. Therefore we may assume $z \in D_2 \setminus \J(x,M)$. By (\ref{ratio}) we have $d(y_i) > |x-y_i|/(e^M-1)$ for $i \in \{ 1,2 \}$ and $d(z) < |x-z|/(e^M-1)$. Since $M > \log 2$ we have $c = 1/(e^M-1) < 1$. Now the boundary $\partial G$ is outside $B^n(y_1,c|x-y_1|) \cup B^n(y_2,c|x-y_2|)$ and has a point in $B^n(z,c|x-z|)$, see Figure \ref{fig1}.

\begin{figure}[htp]
    \begin{center}
      \includegraphics[width=6cm]{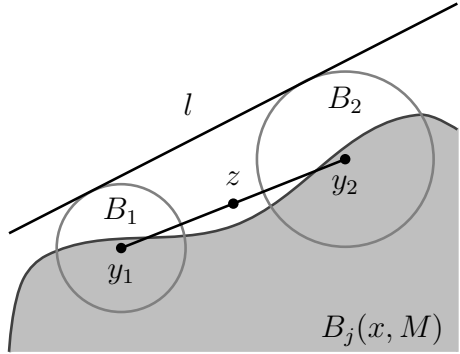}
    \end{center}
    \caption{\label{fig1} Line $l$, Euclidean balls $B_1 = B^n(y_1,c|x-y_1|)$ and $B_2 = B^n(y_2,c|x-y_2|)$ and points $y_1$, $y_2$  and $z$.}
  \end{figure}

  We will show that for $c < 1$ the domain $G$ is not convex. Let us denote by $l$ a line that is a tangent to balls $B^n(y_1,c|x-y_1|)$ and $B^n(y_2,c|x-y_2|)$. Because $d(y_i,l) = c|x-y_i|$ for $i \in \{ 0,1 \}$ we have
  \begin{equation}\label{dzl}
    d(z,l) = \frac{c|x-y_1|+c|x-y_2|}{2}.
  \end{equation}
  By the triangle inequality
  $$
    |x-z| = \left| \frac{x-y_1}{2}+\frac{x-y_2}{2} \right| \le \frac{|x-y_1|}{2}+\frac{|x-y_2|}{2}
  $$
  and by (\ref{dzl})
  $$
    d(z,l) = \frac{c}{2}(|x-y_1|+|x-y_2|) \ge c|x-z|.
  $$
  Now the domain $G$ is not convex, which is a contradiction, and each chord with end points in $\J(x,M) \cap D_2$ is contained in $\J(x,M)$.

  Since each chord with end points in $\J(x,M) \cap D_2$ is contained in $\J(x,M)$, $\J(x,M) \cap D_2 \subset B^n(x,r)$, $D_1$ is convex and $\partial \J(x,M) \cap D_1 \subset S^{n-1}(x,r)$ the $j$-ball $\J(x,M)$ is convex.
\end{proof}

\begin{theorem}\label{jstarlikeness}
  Let $M > 0$ and $G \subsetneq \Rn$ be a starlike domain with respect to $x \in G$. Then the $j$-balls $\J(x,M)$ are starlike with respect to $x$.
\end{theorem}
\begin{proof}
  By Theorem \ref{convexj} we need to consider $M > \log (\sqrt{2}+1)$ which is equivalent to $e^M-1 > \sqrt{2}$. Let us divide $G$ into two parts $D_1 = \{ z \in G \colon d(z) \ge d(x) \}$ and $D_2= G \setminus D_1$.

  Similarly as in the proof of Theorem \ref{convexdomj} the boundary $\partial \J(x,M) \cap D_1$ is a subset of a sphere $S^{n-1}(x,r)$ and $\J(x,M) \subset S^{n-1}(x,r)$. Therefore it is sufficient to show that for each $y \in \J(x,M) \cap D_2$ the line segment $[x,y]$ is in $\J(x,M)$.

  We will show that all chords $[x,y]$ for $y \in \J(x,M) \cap D_2$ are contained in $\J(x,M)$. Let us assume, on the contrary, that there exists points $y_1,y_2 \in \big( \partial \J(x,M) \big) \cap D_2$ with $y_1 \in (x,y_2)$ and $z=(y_1+y_2)/2 \notin \overline \J(x,M)$. Let us first assume $z \in D_1$. Now $j_G(x,z) > j_G(x,y_2)$ is equivalent to $|x-z| /d(x) > |x-y_2| / d(y_2)$. By the selection of $y_1$ and $y_2$ we have $|x-z| < |x-y_2|$ and $d(x) > d(y_2)$ implying $|x-z| /d(x) < |x-y_2| / d(y_2)$, which is a contradiction.

  Let us then assume $z \in D_2$. Now
  $$
    \frac{|x-y_1|}{d(y_1)} = \frac{|x-y_2|}{d(y_2)} = e^M-1 < \frac{|x-z|}{d(z)}
  $$

  \begin{figure}[htp]
    \begin{center}
      \includegraphics[width=6cm]{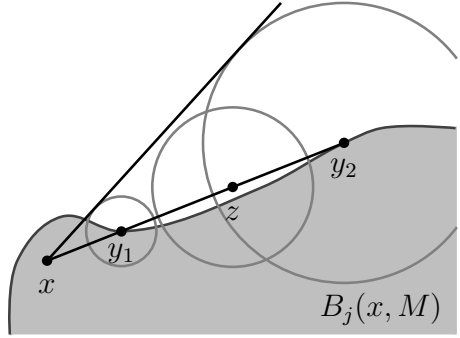}
    \end{center}
    \caption{\label{fig2}Selection of points $y_1$ and $y_2$. Gray circles are $B^n \big( y_1,d(y_1) \big)$, $B^n \big( z,d(z) \big)$ and $B^n \big( y_2,d(y_2) \big)$.}
  \end{figure}
  \noindent and therefore the boundary $\partial G$ does not intersect $B^n \big( y_1,d(y_1) \big)$ or $B^n \big( y_2,d(y_2) \big)$ and contains a point in $B^n \big( z,d(z) \big)$, see Figure \ref{fig2}. This means that $G$ is not starlike with respect to $x$, which is a contradiction.
\end{proof}

\begin{remark}
  (1) Let us consider the domain $G = B^n(0,1) \cup B^n(e_1,1/4) \cup B^n(2e_1,1)$ and show that the $j$-ball $B = B_j(0,\log 3)$ is connected but the $j$-sphere $S = \{ z \in G \colon j_G(0,z) = \log 3 \}$ is disconnected. We have
  $$
    j_G(0,e_1) = \log \left( 1+\frac{1}{1/4} \right) = \log 5
  $$
  and therefore all points $x \in G$ with $x_1 = 1$ are neither in $B$ nor on the boundary $\partial B$. We also have $B, \partial B \subset B^n(0,1) \cup B^n(2e_1,1)$. For all $y \in B^n(2e_1,1) \setminus \{ u \in G \colon \angle{0 \, 2e_1 \, u} < \atan (1/4) \}$ we have
  $$
    j_G(0,y) = \log \left( 1+\frac{|y|}{1-|2-y|} \right) \ge \log \left( 1+2 \right) = \log 3,
  $$
  because $|y|+2|2-y| \ge 2$. For all $y \in B^n(2e_1,1) \cap \{ u \in G \colon \angle{0 \, 2e_1 \, u} < \atan (1/4) \}$ we have
  $$
    j_G(0,y) = \log \left( 1+\frac{|y|}{d(y)} \right) \ge \log \left( 1+\frac{|y_1|}{d(y_1)} \right) \ge \log \left( 1+2 \right) = \log 3
  $$
  and therefore $B \subset B^n(0,1)$ and it is connected.

  Let us now consider $S$ and denote $z \in S$. If $z \in B^n(2e_1,1)$, then $z = 2e_1$. If $z \in B^n(0,1)$, then $z \in \partial B$. Now $S = \partial B \cup \{ 2e_1 \}$ and it is disconnected. In particular, we see that
  $$
    \overline{ \{ z \in G \colon j_G(0,z) < \log 3 \} } \ne \{ z \in G \colon j_G(0,z) \le \log 3 \}.
  $$

  (2) We have seen that in convex domains the $j$-balls are convex and in starlike domains the $j$-balls are starlike. However in simply connected domains the $j$-balls need not be simply connected. Let us consider $G = B^n(0,1) \cup B^n(e_1,h) \cup B^n(2e_1,1)$ for $h \in (0,1)$. Clearly $G$ is simply connected. Let us consider $B = \J(0,\log 4)$. We have
  $$
    j_G(0,2e_1) = \log \left( 1+\frac{2}{1} \right) = \log 3
  $$
  and therefore $2e_1 \in B$. Let $x = (x_1, \dots ,x_n) \in G$ with $x_1 = 1$. Now
  $$
    j_G(0,x) \ge j_G(0,e_1) = \log \left( 1+\frac{1}{h} \right)
  $$
  and $x \notin B$ for $h < 1/3$. For $h = 1/4$ the $j$-ball $B$ is not even connected. Instead of the radius $\log 4$ we could choose any $r > \log 3$.
\end{remark}

\begin{questions}
  We pose some open questions concerning the quasihyperbolic metric and quasihyperbolic balls.
  \begin{itemize}
    \item[(1)] Is it true that for any domain $G \subsetneq \Rn$ and $x \in G$ the quasihyperbolic ball $D_G(x,M)$ is strictly convex if $M \in (0,1]$?
    \item[(2)] Is it true that for any domain $G \subsetneq \Rn$ and $x \in G$ the quasihyperbolic ball $D_G(x,M)$ is strictly starlike with respect to $x$ if $M \in (0,\kappa]$ for $\kappa \approx 2.83297$?
    \item[(3)] Are the quasihyperbolic geodesics unique in every simply connected domain $G \subsetneq \R^2$?
  \end{itemize}
  For the case $\PS$ see Remarks \ref{rem1} and \ref{rem2}.
\end{questions}

\emph{Acknowledgements.} This paper is part of the author's PhD thesis, currently written under the supervision of Prof. M. Vuorinen and supported by the Academy of Finland project 8107317.

%%%%%%%%%%%%%%%%%%%%%%%%%%%%%%%%%%%%%%%%%%%%%%%%%%%%%%%%%%%%%

{\noindent Department of Mathematics\\
University of Turku\\
FIN-20014\\
FINLAND\\
e-mail: riku.klen@utu.fi}

\end{document}